\documentclass[Threeside,11pt]{article}
\usepackage{etoolbox}
\makeatletter
\patchcmd{\thebibliography}
{\list}
{\small          
	\setstretch{1.2}
	\list}
{}{}
\makeatother
\usepackage{setspace}
\usepackage{tikz}
\usepackage{float}
\definecolor{myPurple}{RGB}{109,90,207}
\definecolor{myGreen}{RGB}{4,130,67}
\definecolor{myGold}{RGB}{253,177,71}
\definecolor{myBurgundy}{RGB}{63,1,44}
\definecolor{myTeal}{RGB}{8,119,127}
\definecolor{myCream}{RGB}{255,216,177}
\definecolor{myOrange}{RGB}{225,119,1}

\usetikzlibrary{arrows.meta,bending,calc}
\tikzset{
	ball3Donly/.style={
		circle, minimum size=1.4cm,
		shading=ball,
		ball color=myPurple,
		draw=myPurple
	}
}
\pgfmathsetmacro{\wNear}{2.0pt}
\pgfmathsetmacro{\wMid}{1.3pt}
\pgfmathsetmacro{\wFar}{0.7pt}
\pgfmathsetmacro{\bend}{6mm}
\usepackage [latin1]{inputenc}
\usepackage[english]{babel}
\usepackage{amsmath}
\usepackage{amsthm}
\usepackage{algorithm}
\usepackage{cite}
\usepackage{amsfonts}
\usepackage{amssymb}
\usepackage{mathrsfs}
\usepackage{graphicx}
\usepackage{colortbl,dcolumn}
\usepackage{marvosym}
\usepackage{ifsym}
\usepackage{xcolor}
\usepackage{enumitem}
\allowdisplaybreaks
\usepackage[
colorlinks,
citecolor=blue,
linkcolor=blue]{hyperref} 
\newtheorem{theorem}{\bf Theorem}[section]

\newtheorem{definition}{\bf Definition}[section]
\newtheorem{remark}{\bf Remark}[section]

\numberwithin{equation}{section}

\definecolor{mydeepgreen}{RGB}{3,100,50}
\usepackage{enumitem}
\begin{document}
\title{{\sl Sharp regularity of a weighted  Sobolev space over $ \mathbb{T}^n $ and its relation to finitely differentiable KAM theory}}
\author{Zhicheng Tong\thanks{School of Mathematics, Jilin University, Changchun 130012, P. R.  China. Email: \url{tongzc25@jlu.edu.cn}} 
	\and 
	Yong Li\thanks{Corresponding author. School of Mathematics, Jilin University, Changchun 130012, P. R.  China;   Center for Mathematics and Interdisciplinary Sciences,  Northeast  Normal  University, Changchun 130024, P. R. China. Email: \url{liyong@jlu.edu.cn}}
}

\date{}
\maketitle

\begin{abstract}
In this paper, we investigate the sharp regularity properties of a special weighted Sobolev space defined on the $ n $-dimensional torus, which is of independent interest.  As a key application, we show that for almost all $ n $-dimensional vector fields, the Kolmogorov-Arnold-Moser (KAM) theory holds via this regularity, and in this case, the perturbation must have classical derivatives up to order $ \left[ {  n/2} \right] $, yet it can admit unbounded weak derivatives from order $ \left[ {n/2} \right]+1 $ to $ n$. This result may appear surprising within the classical framework of KAM theory. We also provide further discussion of historical KAM theorems and relevant counterexamples. These findings constitute a new step in the long-standing KAM regularity conjecture.
 \\
	\\
	{\bf Keywords:} {Weighted Sobolev space on the torus,  classical derivatives,  unbounded weak derivatives,   KAM theory with sharp regularity\vspace{2mm}}
\end{abstract}

\tableofcontents

\section{Introduction}
\setcounter{footnote}{0}
\renewcommand{\thefootnote}{{\arabic{footnote}}}
The celebrated Kolmogorov-Arnold-Moser (KAM) theory \cite{MR0068687,MR0163025,MR0147741}, initiated over 70 years ago, primarily concerns the persistence of invariant tori within nearly integrable Hamiltonian systems. Since its inception, KAM theory has played a significant role in various fields, including celestial mechanics, and has evolved into a well-developed array of variants to address a multitude of fundamental problems in dynamical systems and PDEs. Determining the minimal regularity required for KAM theory is indeed a long-standing problem, originating with Moser \cite{MR0147741}. Despite substantial progress over the decades, numerous open questions remain. \textit{Building on these contributions, the present paper establishes a result that may appear counterintuitive.}

Given that the vast majority of KAM theorems rely on analyticity and Newton-type rapid iteration schemes (with super-exponential convergence), the quasi-periodic or almost periodic solutions obtained typically possess a high degree of regularity, such as analyticity (perhaps with some loss in the radius of analyticity). Researchers in KAM theory generally do not focus on \textit{weaker forms of regularity}, such as (weighted) Sobolev regularity\footnote{Although Herman remarked at the 1998 ICM   that ``For many reasons it is not unnatural to work in the Sobolev $ W^{k,p} $ topology'' when considering invariant tori \cite{Her98}, most work in KAM theory has yet to engage with such weak regularity.}, since they often render the KAM iteration ineffective. For instance, without uniform convergence after  infinitely many iterations, one cannot obtain the desired KAM conjugacy directly. It appears that the only known result in PDEs is the recent work by Biasco \textit{et al.} \cite{MR4666300}\footnote{We note that the claim regarding minimal regularity for KAM theory there is slightly inaccurate in view of \cite{POSCHELARXIV}.}, which deals with  weak almost periodic solutions admitting only Sobolev regularity both in time and space for a  one-dimensional nonlinear Schr\"odinger equation, although this is significantly different from the classical KAM theory of Hamiltonian systems.

 When considering the sharp regularity that KAM theory allows for finite-dimensional vector fields, P\"oschel \cite{POSCHELARXIV} constructed a special weighted Sobolev space $ \widetilde {\mathcal{L}_{\tau  + 1,b}^2}({\mathbb{T}^n}) $ (see \eqref{PTSOB} in Section \ref{S1}) on the standard $ n $-dimensional torus $ \mathbb{T}^n: = {\mathbb{R}^n}/ 2\pi {\mathbb{Z}^n} $  in an insightful manner,  and  achieved KAM conjugacy unexpectedly, which fills the gap between the best known KAM results and historical KAM counterexamples\footnote{Although this work remains a preprint (possibly never submitted for publication), this neither implies a flaw in the proof nor diminishes the significance of the results!}. For further details, one can refer to Section \ref{PTSEC2}. It is intriguing to note that the regularity inherent in this weighted Sobolev space is markedly different from classical pointwise differentiability. Not only does it reduce $ 1 $-order classical differentiability required for the best known KAM theory in a weak, non-traditional sense (which we will detail in Remark \ref{PTRE} below Theorem \ref{THPOSCHEL}), as claimed by P\"oschel, \textit{but it also has the capacity to reduce the order of classical differentiability to approximately $ n/2 $, which constitutes one of the main contributions presented in this paper.}  This provides new insights into the \textit{long-standing conjecture} regarding the minimal regularity required in KAM theory.  Moreover, this paper is also driven by the motivation to illustrate the interplay between this special regularity and both  historical KAM theorems and counterexamples (converse KAM theorems). Such illustrations are particularly necessary given that researchers in KAM theory may not be fully acquainted with the work of those who specialize in the destruction of tori. 

\subsection{Sharp regularity of the weighted Sobolev space $ \widetilde {\mathcal{L}_{\tau  + 1,b}^2}({\mathbb{T}^n}) $}\label{S1}
Before proceeding, we shall first introduce in this section the definition of the aforementioned weighted Sobolev space $ \widetilde {\mathcal{L}_{\tau  + 1,b}^2}({\mathbb{T}^n}) $, and establish several sharp results regarding its regularity properties. It is important to emphasize that these results hold independent interest, and the associated technologies may provide precise explanations for some potential results in other fields.

Throughout this paper, we assume that $2 \leqslant b \in \mathbb{N}^+ $ and $ \tau \geqslant n-1 $ with $ 2 \leqslant n \in \mathbb{N}^+ $. Define the weighted Sobolev space $ \widetilde {\mathcal{L}_{\tau  + 1,b}^2}({\mathbb{T}^n}) $ as
\begin{equation}\label{PTSOB}
	\widetilde {\mathcal{L}_{\tau  + 1,b}^2}({\mathbb{T}^n}): =  \left\{ {f: \quad \|f\|_{\widetilde {\mathcal{L}_{\tau  + 1,b}^2}({\mathbb{T}^n})}:=\sum\limits_{\nu  = 0}^\infty  {{{\left( {\sum\limits_{{b^{\nu  - 1}} < \left| k \right| \leqslant {b^\nu }} {{{\left| {{f_k}} \right|}^2}{{\left| k \right|}^{2\tau+2}}} } \right)}^{\frac{1}{2}}}}  <  + \infty } \right\},
\end{equation}
where $ f\left( x \right) = \sum\nolimits_{k \in {\mathbb{Z}^n}} {{f_k}{{\rm e}^{{\rm i}\left\langle {k,x}  \right\rangle}}} \in \mathbb{R}^n $ is a real-valued mapping defined on the standard $ n $-dimensional torus $ \mathbb{T}^n $ with Fourier expansion,  the inner product $ \left\langle {k,x} \right\rangle : = \sum\nolimits_{j = 1}^n {{k_j}{x_j}}  $ and the Fourier coefficients $ {f_k}: = \int_{{\mathbb{T}^n}} {f\left( x \right){{\rm e}^{ - {\rm i}\left\langle {k,x}  \right\rangle}}{\rm d}x}  $, and $ \left| k \right|: = \sum\nolimits_{j = 1}^n {\left| {{k_j}} \right|}  $. 

It is natural that one should consider employing more direct weighted Sobolev spaces  to investigate the regularity of $ \widetilde {\mathcal{L}_{\tau  + 1,b}^2}({\mathbb{T}^n}) $. To this end, for a non-decreasing function $ \varphi(x) \geqslant 0$, we define the weighted Sobolev space $ \mathcal{H}_\varphi ^{\tau+1}({\mathbb{T}^n}) $ as
\[	\mathcal{H}_\varphi ^{\tau+1}({\mathbb{T}^n}): = \left\{ {f: \quad \|	f\|_{\mathcal{H}_\varphi ^{\tau+1}({\mathbb{T}^n})}:=\left(\sum\limits_{k \in {\mathbb{Z}^n}} {{{\left| {{f_k}} \right|}^2}{{\left| k \right|}^{2\tau+2}}\varphi \left( {\left| k \right|} \right)}  \right)^{\frac{1}{2}}<  + \infty } \right\}.\]
In particular, $ \mathcal{H}_\varphi ^{\tau+1}({\mathbb{T}^n}) $ corresponds to the well-known Sobolev space $  \mathcal{H} ^{\tau+1}({\mathbb{T}^n}) $ when $ \varphi(x) $ is a positive constant.   Below, we establish an abstract theorem, which provides a sufficient condition to construct subspaces of $ \widetilde {\mathcal{L}_{\tau  + 1,b}^2}({\mathbb{T}^n}) $.

\begin{theorem}\label{JZT1}
Assume 
	\[\sum\limits_{\nu  = 0}^\infty  {\frac{1}{{\varphi \left( {{b^{\nu  - 1}}} \right)}}}  <  + \infty .\]
Then
\[\mathcal{H}_\varphi ^{\tau+1}({\mathbb{T}^n}) \subset \widetilde {\mathcal{L}_{\tau  + 1,b}^2}({\mathbb{T}^n}).\]
\end{theorem}

In addition to the classical Sobolev spaces $ \mathcal{H}^s (\mathbb{T}^n)$  for $ s >0 $  defined on the torus $ \mathbb{T}^n $ (see B\'{e}nyi and Oh  \cite{MR3119672} for instance), there also exist modified weighted Sobolev spaces incorporating logarithmic terms. For $ \alpha>0 $, define 
 \[\mathcal{H}_{{{\log }^\alpha }}^{\tau+1}({\mathbb{T}^n}): = \left\{ {f: \quad \|f\|_{\mathcal{H}_{{{\log }^\alpha }}^{\tau+1}({\mathbb{T}^n})}:=\left(\sum\limits_{k \in {\mathbb{Z}^n}} {{{\left| {{f_k}} \right|}^2}{{\left| k \right|}^{2\tau+2}}{{\left( {\ln \left( {1 + \left| k \right|} \right)} \right)}^\alpha }}\right)^{\frac{1}{2}}  <  + \infty } \right\}.\]
Based on Theorem \ref{JZT1}, we  are able to deduce a sharp result concerning  the relation between  $ \mathcal{H}_{{{\log }^\alpha }}^{\tau+1}({\mathbb{T}^n}) $ and  $ \widetilde {\mathcal{L}_{\tau  + 1,b}^2}({\mathbb{T}^n}) $.
\begin{theorem}\label{JZT2}
	We have 
	\[   \mathcal{H}_{{{\log }^\alpha }}^{\tau+1}({\mathbb{T}^n})\subset\widetilde {\mathcal{L}_{\tau  + 1,b}^2}({\mathbb{T}^n})\]
	if and only if $ \alpha>1 $.
\end{theorem}

Considering the significance of the Sobolev space $ \widetilde {\mathcal{L}_{\tau  + 1,b}^2}({\mathbb{T}^n}) $ in KAM theory as previously mentioned, \textit{we are deeply interested in the highest order of classical derivatives that mappings within  $ \widetilde {\mathcal{L}_{\tau  + 1,b}^2}({\mathbb{T}^n}) $ possess.} Denote by $ [x] $ the integral part of a real number $ x $.  Theorem \ref{JZT3} presented below provides an affirmative answer.

\begin{theorem}\label{JZT3}
It holds 
\[ \widetilde {\mathcal{L}_{\tau  + 1,b}^2}({\mathbb{T}^n}) \subset {{C}^{\left[ {\tau  + 1 - n/2} \right]}}({\mathbb{T}^n}),\]
where $ C^\ell({\mathbb{T}^n})  $ denotes the standard H\"older space for $ \ell\geqslant 0 $.
\end{theorem}

Interestingly, we observe that Theorem \ref{JZT3} is actually sharp. In other words, \textit{it is indeed possible to construct specific examples from $ \mathbb{T}^n $ to $ \mathbb{T}^n $ where their higher order derivatives are merely weak}\footnote{Throughout this paper, the word ``weak'' means in the $ L^2 $ sense on the torus $ \mathbb{T}^n $. On this aspect, one can refer to 	B\'{e}nyi and Oh \cite{MR3119672} for more details.} \textit{derivatives---they can be essentially  unbounded.} As a sharp result, our  Main Theorem I on the regularity of $  \widetilde {\mathcal{L}_{\tau  + 1,b}^2}({\mathbb{T}^n}) $ in this paper is stated as follows, which has independent interest:
\begin{theorem}[\textbf{Main Theorem I}]\label{JZT4}
Any element in $ \widetilde {\mathcal{L}_{\tau  + 1,b}^2}({\mathbb{T}^n}) $ must possess classical derivatives up to order  $ \left[ {\tau   - n/2} \right] + 1$; however, it can admit unbounded weak derivatives\footnote{The word ``can'' here  indicates sufficiency rather than necessity: there indeed exist elements whose weak derivatives of order $ \left[ {\tau   - n/2} \right] + 2$ up to order $  [\tau]+1 $   are unbounded. Although, as explained earlier, ``weak'' here simply means $ L^{2} $-integrability, we shall present even stronger singularity in Section \ref{SEC22} (Theorem \ref{sssss}).} from order $ \left[ {\tau   - n/2} \right] + 2$ to $  [\tau]+1 $. Moreover, there exist elements in $ \widetilde {\mathcal{L}_{\tau  + 1,b}^2}({\mathbb{T}^n}) $ that do not admit the $ ([\tau]+2) $-th order weak derivatives.
\end{theorem}
\begin{remark}
Consequently,  the highest  classical differentiable order $ \left[ {\tau   - n/2} \right] + 1$ and  the highest weak differentiable order $ [\tau]+1$ in $ \widetilde {\mathcal{L}_{\tau  + 1,b}^2}({\mathbb{T}^n}) $ are both sharp.
\end{remark}
\begin{remark}
It was claimed in \cite{POSCHELARXIV} that for any $ b\geqslant2 $ and $ 2 \leqslant n \in \mathbb{N}^+ $, the $ (n -1) $-th order derivatives of some element in $  \widetilde {\mathcal{L}_{n,b}^2}({\mathbb{T}^n}) $ (i.e., $ \tau=n-1 $) need not be continuous. However, it turns out that when $ n=2 $, this is not the case. 
\end{remark}

Although the decay of Fourier coefficients can be more clearly discerned in the weighted Sobolev space $ \widetilde {\mathcal{L}_{\tau  + 1,b}^2}({\mathbb{T}^n}) $, as demonstrated by Theorem \ref{JZT4}, translating this into regularity is not particularly straightforward. Therefore, we also aim to investigate regularity beyond the classical differentiable order, such as providing a weaker sense of the modulus of continuity as a \textit{sufficient} characterization for $ \widetilde {\mathcal{L}_{\tau  + 1,b}^2}({\mathbb{T}^n}) $ with $ \tau $ being an integer (one can also similarly discuss the case where $ \tau \notin \mathbb{N}^+ $, see   \cite{MR290095} for instance).

Unlike the classical pointwise modulus of continuity (see for instance,  \cite{MR0538680,MR1036903,MR2350326,TLCCM}), we introduce in a weaker sense the modified $ L^2 $-modulus of continuity on the torus $ \mathbb{T}^n $, which can describe the singularity (irregularity) of unbounded functions in the integral sense (also weaker than the classical $ L^2 $-modulus of continuity \cite[Chapter V]{MR290095}). 
\begin{definition}[Modified $ L^2 $-modulus of continuity]
	Given a mapping $ f(x) \in {L^2}\left( {{\mathbb{T}^n}} \right) $,  its $ 0 $-th order modified $ L^2 $-modulus of continuity is defined as 
	\[\widetilde\varpi \left( {f,x} \right): = \mathop {\sup }\limits_{\left| s \right| \leqslant x} {\left( {\int_{{\mathbb{T}^n}} {{{\left| {f\left( {y + s} \right) + f\left( {y - s} \right) - 2f\left( y \right)} \right|}^2}{\rm d}y} } \right)^{\frac{1}{2}}}.\]
	If $ f (x)\in \mathcal{H}^m (\mathbb{T}^n) $ with $ m \in \mathbb{N}^+ $,  its $ m $-th order modified $ L^2 $-modulus of continuity is defined as 
	\[{\widetilde\varpi _{{L^2}}}\left( {{D^m}f,x} \right): = {\left( {\sum\limits_{\left| \alpha  \right| = m} {\widetilde\varpi _{{L^2}}^2\left( {{\partial ^\alpha }f,x} \right)} } \right)^{\frac{1}{2}}},\]
	where $ \alpha  = \left( {{\alpha _1},{\alpha _2},\ldots,{\alpha _n}} \right) $ is a non-negative integral vector, and $ \left| \alpha  \right|: = \sum\nolimits_{j = 1}^n { {{\alpha _j}} } $.
\end{definition}

\begin{definition}
	We say that a modified $ L^2 $-modulus of continuity $ \widetilde\varpi _{{L^2}}\left( {  \cdot  ,x} \right) $ is of Dini type if 
	\[\int_0^1 {\frac{{{{\widetilde \varpi }_{{L^2}}}\left( { \cdot ,x} \right)}}{x}{\rm d}x}  <  + \infty. \]
\end{definition}

\begin{theorem}\label{TPDLINGLI4}
	If $ 	f(x): \mathbb{T}^n \to \mathbb{R}^n $ admits a $ (\tau+1)  $-th order Dini modified $ L^2 $-modulus of continuity, where $ 2 \leqslant b \in \mathbb{N}^+ $ and $ 1 \leqslant n-1 \leqslant \tau \in \mathbb{N}^+ $,  then $ 	f(x)\in \widetilde {\mathcal{L}_{\tau  + 1,b}^2}({\mathbb{T}^n})  $. 
\end{theorem}

Particularly, one might consider integrating the techniques from  Zygmund  \cite{MR1963498} and Fan and Meyer \cite{MR4529026}---of course, with the non-trivial task of generalizing these methods from $ \mathbb{T}^1 $ to $  \mathbb{T}^n $---to estimate the asymptotic behavior  of elements within $ \widetilde {\mathcal{L}_{\tau  + 1,b}^2}({\mathbb{T}^n}) $  at singularity points based on the decay of their Fourier coefficients;  however, we will not delve into this exploration here. In Section \ref{SEC22}, we further explore  the possible singularity of the weak derivatives of mappings within  $\widetilde{\mathcal{L}_{\tau + 1,b}^2}(\mathbb{T}^n)$, as detailed in Theorem \ref{sssss}; however, we prefer  not to state the result here to prevent excessive verbosity.


\subsection{Application to the finitely differentiable KAM theorem}\label{PTSEC2}

This section is devoted to providing a review of historical KAM theorems and counterexamples, with a particular emphasis on the aspect of regularity, especially within the context of finite differentiability. Additionally, this section includes discussions aimed at illustrating the importance of our Main Theorems (Theorem \ref{JZT4} and \ref{THPOSCHEL}) for the reader's understanding.

We continue to assume that $2 \leqslant b \in \mathbb{N}^+ $ and $ 2 \leqslant n \in \mathbb{N}^+ $.  We begin by presenting the definition of the Diophantine nonresonance  that characterizes the arithmetic properties of frequencies. This concept is fundamental in KAM theory.

\begin{definition}
We say that a vector $ \omega \in \mathbb{R}^n $ is Diophantine with exponent $ \tau>0 $, if for some $ \alpha>0 $,
\[|\langle k, \omega\rangle| \geqslant \frac{\alpha}{|k|^{\tau}}, \quad \forall 0 \ne k \in \mathbb{Z}^n.\]
\end{definition}
\begin{remark}
If the vector $ \omega \in \mathbb{R}^n  $ does not satisfy the Diophantine condition, we classify it as Liouville.
\end{remark}

It is well known (see R\"ussmann \cite{MR467824} for instance) that when $ \tau > n-1 $, such vectors have full Lebesgue measure in $ \mathbb{R}^n $; when $ \tau = n-1 $, they have zero Lebesgue measure but possess Hausdorff dimension $ n $, hence are continuum many; and when $ \tau < n-1 $, they do not exist according to Dirichlet's principle. 

It is well known that the two essential assumptions in KAM theory are the regularity of the system and the nonresonance of frequencies, both of which are indispensable; otherwise, the invariant tori may break. Over the past 70 years since the birth of KAM theory, many efforts have been made to weaken the regularity of the system.

Notably, Moser's pioneering work \cite{MR0147741} reduced the analyticity utilized by Kolmogorov and Arnold \cite{MR0068687,MR0163025} to finite differentiability of $ C^{333} $, focusing on area-preserving mappings with the monotone twist property (this can also  correspond to a non-analytic Hamiltonian version with  $ 2 $-degrees of freedom). 	Herman \cite{MR0728564} constructed a $ C^{3-\varepsilon} $ counterexample to show 	the nonexistence of an invariant curve, where $ \varepsilon>0 $ is arbitrarily given (here and below). Later he provided a general   counterexample in  \cite{MR0874026} with $ C^{n+1-\varepsilon} $ regularity in $ n $ dimensions.  We also mention the counterexamples with lower regularity constructed by Takens \cite{MR300311} and	Mather \cite{MR0766107}.		Following Moser's approximating idea \cite{MR0265692} and  P\"{o}schel's work \cite{MR0595866},	Salamon \cite{MR2111297} considered the persistence of an individual torus. He proved that, for a given non-degenerate Hamiltonian system with $ n \geqslant2  $-degrees of freedom  and a Diophantine frequency with exponent $ \tau>n-1 $, both the unperturbed Hamiltonian (not necessarily integrable) and the perturbation can be of class $ C^{\ell} $ with $ \ell>2\tau+2 $ (therefore $ \ell>2n $), and the frequency-preserving torus is $ C^{s+\tau} $, the dynamic on it is $ C^{s} $ conjugated to the linear flow, where $ s \leqslant m+1 $ with $ m <\ell-2\tau-2 $, and $ s, s+\tau \notin \mathbb{N} $. This is a nearly sharp result, considering that 	Cheng and Wang \cite{MR3061774} constructed a   counterexample: for an integrable Hamiltonian with  $ n \geqslant 2 $-degrees of freedom, any Lagrangian torus with a given unique rotation vector can be destroyed by arbitrarily $ C^{2n-\varepsilon} $-small perturbations.  	Albrecht \cite{MR2350326} also  investigated a non-degenerate Hamiltonian system, but for a non-universal Diophantine frequency with exponent $ \tau=n-1 $ (with zero Lebesgue measure). He proved that $ C^{2n} $ plus a Dini modulus of continuity $ \varpi $ (i.e., $ \int_0^1 {\frac{{\varpi \left( x \right)}}{x}{\rm d}x}  <  + \infty  $) is sufficient to obtain the persistent frequency-preserving torus. 	However, unlike Salamon's work \cite{MR2111297}, the unperturbed part here is analytic and integrable, while the remaining regularity has not been investigated. Independently of Albrecht's work \cite{MR2350326},	the authors \cite{TLCCM} extended Salamon's universal KAM frequency-preserving persistence \cite{MR2111297} in the same setting, weakened the previous H\"older type regularity to  $ C^{[2\tau]+2} $ plus a general Dini type modulus of continuity depending on $ \tau>n-1 $ (in particular, if $ \tau=n-1 $, then it recovers  \cite{MR2350326}). 	Moreover, the authors employed  asymptotic analysis for the first time to investigate the remaining regularity regarding the modulus of continuity, removing some restrictions on parameters (e.g., $ s $, $ m $ and $ \tau $) at the expense of  quantitative estimates of toroidal deformations. 	Furthermore, the authors examined cases with higher differentiability.

	In addition to the aforementioned KAM work focusing on individual tori, there are also studies that consider the KAM persistence for sets with asymptotically full Lebesgue measure. In this regard, the pioneering work originated from Lazutkin \cite{MR0328219}, who investigated twist mappings  with excessively high regularity requirements. Considering the Hamiltonian system with  $ n \geqslant 2 $-degrees of freedom,	P\"{o}schel \cite{MR0668410} reduced the regularity of the perturbation to \( C^{\ell+\tau} \), where \( \ell > 2\tau + 2 \) and \( \tau > n - 1 \), while demanding analyticity for the unperturbed part. Bounemoura \cite{MR4167789} further improved upon P\"{o}schel's result \cite{MR0668410} by requiring the perturbation to be \( C^{\ell} \) (as in Salamon \cite{MR2111297}), and the unperturbed part to be \( C^{\ell+2} \).  Koudjinan \cite{MR4104457} developed Bounemoura's findings \cite{MR4167789}, lowering the regularity of the unperturbed part to \( C^{\ell} \), but at the cost of less favorable measure estimate for the complement of the
set of tori. 	Undoubtedly, building on the ideas from Albrecht \cite{MR2350326} and the authors \cite{TLCCM}, the regularity requirements for KAM persistence of positive measure can be further reduced in terms of the modulus of continuity. However, there is little room left for further progress, at least for the persistence of an individual torus; hence, each step forward might be  challenging.

Based on finite differentiability within finite-dimensional settings, many other  efforts in  KAM theory, as well as destruction theory, have been made in various contexts beyond classical Hamiltonian systems with maximal tori. This involves historical and very recent progress, for instance, see Chierchia and Qian \cite{MR2093919}, Wagener \cite{Wag10}, Khesin \textit{et al.} \cite{MR3269186}, Wang \cite{MR3111869, MR4385768, WLARXIV}, Huang \textit{et al.} \cite{MR3708154}, Li and Shang \cite{MR3960504}, Li \textit{et al.} \cite{MR4618673}, Hu \cite{MR4623276},   Hu and Zhang \cite{MR4714539}, Sorrentino and Wang \cite{SW26} and the references therein. While in infinite-dimensional settings, the  authors \cite{arXiv:2306.08211} have for the first time established a $ C^\infty $ type KAM theorem. This is indeed sharp, as  finite differentiability would destroy the KAM persistence in this case, as previously mentioned.

As KAM theory has expanded into a diverse and intricate domain, efforts to construct an all-encompassing framework that addresses every conceivable scenario have been made, yet a gap persists. This gap is particularly pronounced in the context of finite differentiability, where standard methods frequently fall short of the required precision. For instance, via Diophantine frequencies with exponent $\tau > n-1$, the regularity needed for the finite differentiability case, as discussed in the reference \cite{ASARXIV} by Alazard and Shao, exceeds the sharp differentiability (see  \cite{MR2111297,TLCCM}) by approximately $n/2$ orders, with $ n $ representing the number of degrees of freedom in the Hamiltonian system. Consequently, the pursuit of novel methodologies to delve into KAM theory with finite differentiability remains of great importance and is not readily supplanted by standard approaches.

Back to our concern on the  minimum regularity necessary  for KAM theory, 	P\"{o}schel  commented in \cite{POSCHELARXIV}, \textit{``Consider the case of $ n $-degrees of freedom and $ \tau=n-1 $. It was conjectured for a long time, and sometimes even stated as fact, that $ C^n $ is the minimal regularity requirement for KAM to apply in this case.''} On this aspect, one can refer to  \cite{MR0874026,MR2350326,MR3061774,MR4167789,MR4666300}, for instance. However, as an unexpected breakthrough, he introduced the weighted Sobolev space $ 	\widetilde {\mathcal{L}_{\tau  + 1,b}^2}({\mathbb{T}^n}) $  defined in \eqref{PTSOB} and proved that if the perturbation $ P(x) $ of the  simplest yet fundamental perturbed system 
\begin{equation}\label{PTPO}
	x' = \omega  + P\left( x \right),\quad x \in {\mathbb{T}^n}, \quad \text{$ \omega $ is  Diophantine with exponent $ \tau \geqslant n-1 \geqslant 1$}
\end{equation}
is sufficiently small in the  $ 	\widetilde {\mathcal{L}_{\tau  + 1,b}^2}({\mathbb{T}^n}) $ norm,  then there exists a modifying term $ \widetilde{\omega} $ and a nearly identity transformation $ \Psi $ such that $ \omega -\widetilde \omega + {P } $ is conjugated to $ \omega $:
\[	\Psi^ * \left( {\omega -\widetilde \omega + {P } } \right) = \omega ;\] 
this process is also known as the linearization of the vector field on the torus. Consequently, as an immediate  application, he claimed that Moser's  KAM theory with parameters could apply to small perturbations of weaker regularity than $ C^n $, yet  still allows for Herman's $ C^{n-\varepsilon} $ (for arbitrary $ \varepsilon>0 $) counterexample; see also Section \ref{SEC21}. However, the authors observe that while P\"{o}schel presented the above KAM theorem which is almost unimprovable, he did not elucidate the regularity of the perturbation with sufficient precision. The main motivation behind formulating Main Theorem I is to provide some novel insights into this long-standing KAM regularity conjecture, as detailed in Main Theorem II below.

By combining P\"{o}schel's finitely differentiable KAM theorem and Main Theorem I (Theorem \ref{JZT4}), we arrive at the second main result in this paper:

\begin{theorem}[\textbf{Main Theorem II}]\label{THPOSCHEL}
(I) Assume that the frequency $\omega$ of the vector field in \eqref{PTPO} is Diophantine with exponent $ n-1 $. Then there exist a modifying term $ \widetilde \omega $ and a nearly identity transformation $ \Psi $, such that $ \omega -\widetilde \omega + {P } $ is conjugated to $ \omega $ as follows:
\begin{equation}\label{PTge}
	\Psi^ * \left( {\omega -\widetilde \omega + {P } } \right) = \omega ,
\end{equation}
 provided that $ P(x) \in   \widetilde {\mathcal{L}_{n,b}^2}({\mathbb{T}^n})  $ is sufficiently small in the norm sense. In this case, the perturbation $ P(x) $ must have classical derivatives up to order  $ \left[ {  n/2} \right] $,  but can admit unbounded weak derivatives from order $ \left[ {n/2} \right]+1 $ to $ n$.\vspace{2mm}
 \\
  (II) Furthermore, the aforementioned  differentiability for the KAM perturbation applies to  almost all fixed frequencies $ \omega \in \mathbb{R}^n $.\vspace{2mm}\\
  (III)	If the perturbation $ P(x) $  admits a small $ n $-th order Dini modified $ L^2 $-modulus of continuity, then the KAM conjugacy in \eqref{PTge} is valid.
\end{theorem}
\begin{remark}\label{PTRE}
	Although P\"oschel's original result shows a  difference in the smallness of the perturbation norm for $ \tau = n-1 $ and $ \tau > n-1\ $ (i.e., the smallness of  $ \|P\|_{   \widetilde {\mathcal{L}_{n,b}^2}({\mathbb{T}^n})} $ and $ \|P\|_{   \widetilde {\mathcal{L}_{\tau,b}^2}({\mathbb{T}^n})} $), we point out that there is no difference from the perspective of differentiability. Therefore, we present a universal KAM result. We insist on providing the first part to illustrate that, based on P\"oschel's results, our differentiability cannot be further weakened, even if we seek better irrational vectors (less close to rational vectors).

\end{remark}

To conclude this section, we provide below a table to help readers better understand Main Theorem II (Theorem \ref{THPOSCHEL}), particularly  from a KAM perspective:
{\large
\begin{equation}\notag
	\renewcommand\arraystretch{2}
	\begin{array}{|c|}
		\hline \textbf {Derivatives of the  KAM perturbation $ P(x) \in   \widetilde {\mathcal{L}_{n,b}^2}({\mathbb{T}^n})$: The weakest case} \\
		\hline \underbrace {1,  \quad \cdots\cdots , \quad \left[ \frac{{n}}{2} \right]}_{\rm{classical\;derivatives}}, \quad  \quad  \quad \underbrace {\left[ \frac{{n}}{2} \right] + 1, \quad  \cdots\cdots ,\;\overbrace{n-1}^{\rm the \; critical\; order}, \quad  n }_{\rm{\Huge unbounded\;weak\;derivatives}} \\
		\hline
	\end{array}
\end{equation}}

\noindent We will postpone the explanation of the criticality of the $(n-1)$-th order to Section \ref{PTSEC4}. Therefore, although we do not enhance the KAM technique itself\footnote{After consulting with some experts, it seems that there is little hope for improving it.}, we improve P\"oschel's claim on the classical differentiable order of the perturbation by about \textit{half}.

The rest of this paper is organized as follows: Section \ref{PTSEC4} discusses  some further relationship between the weighted Sobolev space $ \widetilde {\mathcal{L}_{\tau  + 1,b}^2}({\mathbb{T}^n}) $ and both historical KAM theorems and  counterexamples (Section \ref{SEC21}), as well as the possible singularity of unbounded weak derivatives of elements within $ \widetilde {\mathcal{L}_{\tau  + 1,b}^2}({\mathbb{T}^n}) $ (Section \ref{SEC22}); Section \ref{PTSEC5} is dedicated to presenting the proofs of all the results discussed in this paper, including Theorems \ref{JZT1} to \ref{sssss}.

\section{Further discussions}\label{PTSEC4}
Considering  P\"oschel's contributions to the finitely differentiable KAM theory and  additional insights into regularity presented in this paper (as detailed in Theorem \ref{THPOSCHEL}), one might feel that they are somewhat incompatible with historical KAM theorems and counterexamples, however, this is not the case. We will explain in detail in this section, and this is also one of the key motivations behind the composition of this paper. Such discussions may be highly important for delving deeper into the discrepancies between KAM theory and counterexamples, providing valuable insights to eliminate certain cases.

\subsection{Relation to  historical KAM theorems and  counterexamples}\label{SEC21}

\begin{itemize}
	\item[(A)] 
	KAM theory for differential forms and mapping forms can often be transformed into each other, especially in the finitely differentiable setting (using	the technique of generating functions); see Douady \cite{Dou82a,Dou82b},	Kuksin and P\"oschel \cite{KP94} for instance.  Consequently, P\"oschel indicated that his result can also refine Moser's KAM theorem, which means that the regularity required  for the perturbation can be lower than $ C^n $\footnote{It is critical due to Herman's counterexample \cite{MR0874026}; as in the differential form \eqref{PTPO}, $ C^{n-1} $ is the critical case.}: the $ n $-th order derivatives can possess  weighted Sobolev regularity without necessarily being continuous. Therefore, our Theorem \ref{THPOSCHEL}  also significantly weakens  the mapping form of  KAM theory from the perspective of the classical differentiable order of the perturbation.

	\item[(B)] P\"oschel's result \cite{POSCHELARXIV} is \textit{not} in conflict with the H\"older-type or the modulus of continuity type results established by Herman \cite{MR0728564,MR0874026}, Salamon \cite{MR2111297}, Albrecht \cite{MR2350326}, Cheng and Wang \cite{MR3061774}, 	Bounemoura \cite{MR4167789}, Koudjinan \cite{MR4104457}, and the authors \cite{TLCCM}. For instance, for  Hamiltonian systems with $ n $-degrees of freedom, $ C^{\ell} $ H\"older regularity with $ \ell>2\tau+2 $ (critical) \textit{does not} imply the existence of higher order weak derivatives; whereas Theorem \ref{THPOSCHEL} indicates that P\"oschel's $ 	\widetilde {\mathcal{L}_{\tau  + 1,b}^2}({\mathbb{T}^n}) $ regularity,  although allowing the $ (n-1) $-th (critical) order derivatives of the perturbation in \eqref{PTPO} to be non-continuous, but implies the existence of the $ n $-th  order weak derivatives (with a higher order if it can be generalized to Hamiltonian systems). Specifically, \cite{TLCCM} also constructs a  Hamiltonian system where the highest-order derivatives are nowhere H\"older continuous. Consequently, P\"oschel's result, although being a breakthrough, is independent of some known KAM results \cite{MR2111297,MR2350326,MR4167789,MR4104457,TLCCM} and \textit{does not} exhibit an inclusion relationship.

\item[(C)] Although we can regard Theorem \ref{THPOSCHEL} as a special instance  of the KAM theorem (with attracting invariant tori, see Koch \cite{MR4735770} for instance) for a Hamiltonian function  given by 
\[H\left( {x,y} \right) = \left\langle {\omega ,y} \right\rangle  + \left\langle {P\left( x \right),y} \right\rangle , \quad x \in {\mathbb{T}^n}, \quad y \in G \subset {\mathbb{R}^n},\]
 and in this case, the regularity of the Hamiltonian function $ H\left( {x,y} \right) $ is quite weak, entirely depending on that of $ P(x) $ (i.e., $ H\left( {x,y} \right) $ can    admit unbounded weak derivatives from order $ \left[ {n/2} \right]+1 $ to $ n$), this \textit{does not} encompass the situation for \textit{arbitrary} perturbations $ \widetilde{P}(x,y) $.

\item[(D)] Does   P\"oschel's approach \cite{POSCHELARXIV} apply to nearly integrable Hamiltonian systems?  What about frequency-preserving tori under certain nondegeneracy conditions? How does it perform in degenerate settings? And what about its applicability to lower-dimensional tori? These are far from trivial questions.  If feasible, one might seek some room between  KAM theory (Treshch\"{e}v \cite{Tre89}, Li and Yi \cite{MR2003447,LY05}, Salamon \cite{MR2111297},  Sevryuk \cite{Sev06}, Albrecht \cite{MR2350326}, 	Han \textit{et al.} \cite{MR2639543}, Qian \textit{et al.}  \cite{MR4669322} and the authors \cite{TLCCM}, for instance) and  converse KAM theory (Cheng and Wang \cite{MR3061774}, for instance).  Meanwhile, the results in \cite{POSCHELARXIV} and this paper could prove informative for several other lines of work, for instance MacKay and Percival \cite{MP85}, MacKay \cite{Mac89,Mac18}, MacKay \textit{et al.} \cite{MMS89}, Sevryuk \cite{Sev08,Sev17}, Wagener \cite{Wag10},  Meiss \cite{Mei12},  Duignan and Meiss \cite{DM21}, as well as the references therein.

\item[(E)] Consider the perturbed vector field $ \omega +P(x)$ in Theorem \ref{THPOSCHEL}. It is well known that when the rotation set of $ \omega +P(x)$ contains $ \omega $, the modifying term $ \widetilde{\omega} $ in the KAM conjugacy \eqref{PTge} vanishes, i.e., achieving frequency-preserving. However, this is not always the case for arbitrary vector fields close to $ \omega $. As a consequence, for a fixed Diophantine frequency $ \omega \in \mathbb{R}^n $, if we aim to break up \textit{all} tori admitting frequency $ \omega $ through a \textit{uniform} perturbation $ \mathcal{P}(x) $\footnote{This remains an open problem listed at ICM 2018 concerning the sharpness in higher dimensions.}, then the regularity for $ \mathcal{P}(x)  $ must be lower than that in $ 	\widetilde {\mathcal{L}_{\tau  + 1,b}^2}({\mathbb{T}^n}) $. According to Theorem \ref{THPOSCHEL}, the desired perturbation $ \mathcal{P}(x)  $ may possess  classical derivatives of order lower than $   \left[ {  n/2} \right]  $, but may possess certain higher order weak (unbounded) derivatives.

\item[(F)] It is important to observe that almost all historical KAM counterexamples in the finitely differentiable context were constructed within the H\"older norm. As a consequence, motivated by Theorem \ref{THPOSCHEL},  the question of  whether the breakdown of tori can be further explored  under the modulus of continuity, or even in certain weighted Sobolev norms (higher than  that in $ \widetilde {\mathcal{L}_{\tau  + 1,b}^2}({\mathbb{T}^n}) $), remains open.


\end{itemize}

\subsection{The possible singularity of unbounded weak derivatives}\label{SEC22}
Recall that Theorem \ref{TPDLINGLI4} provides a regularity perspective in terms of the modified $ L^2 $-modulus of continuity. While in this section, we aim to explore  the possible unboundedness of higher order weak derivatives of mappings in $ \widetilde {\mathcal{L}_{\tau  + 1,b}^2}({\mathbb{T}^n}) $. As mentioned below Theorem \ref{TPDLINGLI4}, one may develop techniques from Zygmund  \cite{MR1963498} and Fan and Meyer \cite{MR4529026} to analyze the asymptotic behavior of mappings with in $ \widetilde {\mathcal{L}_{\tau  + 1,b}^2}({\mathbb{T}^n}) $ as well as their derivatives via specific (or  critical) Fourier coefficients. However, we have a more intuitive way of looking at this point, although this approach is somewhat rough  (reflected in the power). Let us consider the simplest case where $ \tau =n-1 $;  more general cases are omitted here. We point out in the following Theorem \ref{sssss} that the higher order weak derivatives of elements within $ {\mathcal{L}_{n,b}^2}({\mathbb{T}^n}) $ can exhibit  asymptotic behavior at some singularity point $ {x_0} \in {\mathbb{T}^n} $ as 
\[\frac{1}{{{{\left| {x - {x_0}} \right|}^\mu }{{\left( { - \ln \left| {x - {x_0}} \right|} \right)}^\eta }}},\quad \text{where $ \mu, \eta>0 $ are some suitable constants.}\]

\begin{theorem}\label{sssss}
For any $ {x_0} \in {\mathbb{T}^n} $ and $ \sigma>1/2 $, there exists $ f(x)  \in \widetilde {\mathcal{L}_{n,b}^2}({\mathbb{T}^n}) $, such that as $ x \to {x_0} $, 
\begin{equation}\label{sss}
	\left|{D^{n - j}}f\left( x \right)\right| \simeq \int_{\left| {x - {x_0}} \right|}^{\frac{1}{2}} { \cdots \int_{\left| {{v_3}} \right|}^{\frac{1}{2}} {\int_{\left| {{v_2}} \right|}^{\frac{1}{2}} {\frac{1}{{v_1^{\frac{n}{2}}{{\left( { - \ln {v_1}} \right)}^\sigma}}}{\rm d}{v_1}{\rm d}{v_2} \cdots {\rm d}{v_j}} } } , \quad 1 \leqslant j \leqslant n,
\end{equation}
and 
\[\left| {{D^n}f\left( x \right)} \right| \simeq \frac{1}{{{{\left| {x - {x_0}} \right|}^{\frac{n}{2}}}{{\left( { - \ln \left| {x - {x_0}} \right|} \right)}^\sigma }}}.\]
\end{theorem}

By employing Theorem \ref{sssss}, one can attain a more intuitive comprehension of the regularity of perturbations in Theorem \ref{THPOSCHEL}, specifically from a quantitative perspective on the unbounded singularity.

\section{Proofs of all results}\label{PTSEC5}
For the sake of brevity, in the subsequent analysis, we will use the notation  $ f \gtrsim g $  (or $ f \lesssim g  $) to mean there exists a universal constant $ C > 0 $, such that $ Cf \geqslant g $ (or $ f \leqslant Cg $). If $ f \gtrsim g $ and $ f \lesssim g $ hold simultaneously, we will simply write $ f \simeq g $.
Below, we provide detailed proofs for all results, from Theorem \ref{JZT1} to Theorem \ref{sssss}.

\subsection{Proof of Theorem \ref{JZT1}}
	Let $ f(x)=\sum\nolimits_{k \in {\mathbb{Z}^n}} {{f_k}{{\rm e}^{{\rm i}\left\langle {k,x}  \right\rangle}}} \in \widetilde {\mathcal{L}_{\tau  + 1,b}^2}({\mathbb{T}^n}) $, then it follows that
\begin{align*}
	\sum\limits_{\nu  = 0}^\infty  {{{\left( {\sum\limits_{{b^{\nu  - 1}} < \left| k \right| \leqslant {b^\nu }} {{{\left| {{f_k}} \right|}^2}{{\left| k \right|}^{2\tau+2}}} } \right)}^{\frac{1}{2}}}}  	& \leqslant \sum\limits_{\nu  = 0}^\infty  {\frac{1}{{ {\varphi^{\frac{1}{2}} \left( {{b^{\nu  - 1}}} \right)} }}{{\left( {\sum\limits_{{b^{\nu  - 1}} < \left| k \right| \leqslant {b^\nu }} {{{\left| {{f_k}} \right|}^2}{{\left| k \right|}^{2\tau+2}}\varphi \left( |k| \right)} } \right)}^{\frac{1}{2}}}} \\
	&  \leqslant \sum\limits_{\nu  = 0}^\infty  {\left( {\frac{1}{{\varphi \left( {{b^{\nu  - 1}}} \right)}} + \sum\limits_{{b^{\nu  - 1}} < \left| k \right| \leqslant {b^\nu }} {{{\left| {{f_k}} \right|}^2}{{\left| k \right|}^{2\tau+2}}\varphi \left( |k| \right)} } \right)} \\
	&  = \sum\limits_{\nu  = 0}^\infty  {\frac{1}{{\varphi \left( {{b^{\nu  - 1}}} \right)}}}  + \sum\limits_{k \in {\mathbb{Z}^n}} {{{\left| {{f_k}} \right|}^2}{{\left| k \right|}^{2\tau+2}}\varphi \left( {\left| k \right|} \right)} \\
	& <  + \infty .
\end{align*}
This proves Theorem \ref{JZT1}.

\subsection{Proof of Theorem \ref{JZT2}}
Note that for $ \alpha>1 $, we have
\[\sum\limits_{\nu  = 0}^\infty  {\frac{1}{{\varphi \left( {{b^{\nu  - 1}}} \right)}}}  = \sum\limits_{\nu  = 0}^\infty  {\frac{1}{{{{\ln }^\alpha }\left( {1 + {b^{\nu  - 1}}} \right)}}}  \simeq \sum\limits_{\nu  = 1}^\infty  {\frac{1}{{{\nu ^\alpha }}}}  <  + \infty. \]
Therefore, by applying Theorem \ref{JZT1}, we prove that $  \mathcal{H}_{{{\log }^\alpha }}^{\tau+1}({\mathbb{T}^n})\subset\widetilde {\mathcal{L}_{\tau  + 1,b}^2}({\mathbb{T}^n}) $.

Next, for any $0< \alpha \leqslant 1 $, let us construct an example $ f(x)=\sum\nolimits_{k \in {\mathbb{Z}^n}} {{f_k}{{\rm e}^{{\rm i}\left\langle {k,x}  \right\rangle}}} $  (depending on $ \alpha $), such that $ f(x) \in \mathcal{H}_{{{\log }^\alpha }}^{\tau+1}({\mathbb{T}^n})$, but $ f(x) \notin \widetilde {\mathcal{L}_{\tau  + 1,b}^2}({\mathbb{T}^n}) $. Let $ {f_k} = 0 $ for  $ \left| k \right| \notin {\left\{ {{b^{\nu  - 1}}} \right\}_{\nu  \in \mathbb{N}}} $, while $ \left| {{f_k}} \right| = h\left( \nu  \right) $ for $ \left| k \right| = {b^{\nu  - 1}} $ for some $ \nu \in \mathbb{N} $. In this case, we have 
\begin{align}
	\sum\limits_{k \in {\mathbb{Z}^n}} {{{\left| {{f_k}} \right|}^2}{{\left| k \right|}^{2\tau+2}}{{\left( {\ln \left( {1 + \left| k \right|} \right)} \right)}^\alpha }}  & = \sum\limits_{\nu  = 1}^\infty  {\left( {\sum\limits_{\left| k \right| = {b^{\nu  - 1}}} {{{\left| {{f_k}} \right|}^2}{{\left| k \right|}^{2\tau+2}}{{\left( {\ln \left( {1 + \left| k \right|} \right)} \right)}^\alpha }} } \right)} \notag \\
	& \simeq \sum\limits_{\nu  = 1}^\infty  {\left( {\left( {\sum\limits_{\left| k \right| = {b^{\nu  - 1}}} 1 } \right){h^2}\left( \nu  \right){b^{(2\tau+2)\nu }}{\nu ^\alpha }} \right)} \notag \\
	& \simeq \sum\limits_{\nu  = 1}^\infty  {\mathscr{S}_n\left( {{b^{\nu  - 1}}} \right){h^2}\left( \nu  \right){b^{(2\tau+2)\nu }}{\nu ^\alpha }} \notag \\
	& \simeq \sum\limits_{\nu  = 1}^\infty  {{b^{\left( {n - 1} \right)\nu }}{h^2}\left( \nu  \right){b^{(2\tau+2)\nu }}{\nu ^\alpha }} \notag \\
	& \simeq \sum\limits_{\nu  = 1}^\infty  {{b^{\left( {2\tau+n +1} \right)\nu }}{h^2}\left( \nu  \right){\nu ^\alpha }} ,\notag 
\end{align}
and
\begin{align*}
	\sum\limits_{\nu  = 0}^\infty  {{{\left( {\sum\limits_{{b^{\nu  - 1}} < \left| k \right| \leqslant {b^\nu }} {{{\left| {{f_k}} \right|}^2}{{\left| k \right|}^{2\tau+2}}} } \right)}^{\frac{1}{2}}}}  &\simeq \sum\limits_{\nu  = 1}^\infty  {{{\left( {\sum\limits_{\left| k \right| = {b^{\nu  - 1}}} {{{\left| {{f_k}} \right|}^2}{{\left| k \right|}^{2\tau+2}}} } \right)}^{\frac{1}{2}}}} \\
	& \simeq \sum\limits_{\nu  = 1}^\infty  {{{\left( {\sum\limits_{\left| k \right| = {b^{\nu  - 1}}} {{h^2}\left( \nu  \right){b^{({2\tau+2})\nu }}} } \right)}^{\frac{1}{2}}}} \\
	& \simeq \sum\limits_{\nu  = 1}^\infty  {{{\left( {\mathscr{S}_n\left( {{b^{\nu  - 1}}} \right){h^2}\left( \nu  \right){b^{({2\tau+2})\nu }}} \right)}^{\frac{1}{2}}}} \\
	& \simeq \sum\limits_{\nu  = 1}^\infty  {{{\left( {{b^{\left( {n - 1} \right)\nu }}{h^2}\left( \nu  \right){b^{({2\tau+2})\nu }}} \right)}^{\frac{1}{2}}}} \\
	& \simeq \sum\limits_{\nu  = 1}^\infty  {{b^{\left( {2\tau+n+1} \right)\nu /2}}h\left( \nu  \right)} ,
\end{align*}
where $ \mathscr{S}_n(x) $ represents the surface area of an $ n $-dimensional sphere with radius $ x $ (or the volume of an $ (n-1) $-dimensional sphere). Taking ${b^{\left( {2\tau+n+1} \right)\nu /2}}h\left( \nu  \right) \sim \frac{1}{{\nu \ln \nu }} $ as $ \nu \to +\infty $, we obtain that
\[	\sum\limits_{k \in {\mathbb{Z}^n}} {{{\left| {{f_k}} \right|}^2}{{\left| k \right|}^{2\tau+2}}{{\left( {\ln \left( {1 + \left| k \right|} \right)} \right)}^\alpha }}  \simeq \sum\limits_{\nu  = 2}^\infty  {\frac{1}{{{\nu ^{2 - \alpha }}{{\ln }^2}\nu }}}  \leqslant \sum\limits_{\nu  = 2}^\infty  {\frac{1}{{\nu {{\ln }^2}\nu }}}  <  + \infty ,\]
and
\[	\sum\limits_{\nu  = 0}^\infty  {{{\left( {\sum\limits_{{b^{\nu  - 1}} < \left| k \right| \leqslant {b^\nu }} {{{\left| {{f_k}} \right|}^2}{{\left| k \right|}^{2\tau+2}}} } \right)}^{\frac{1}{2}}}}  \simeq\sum\limits_{\nu  = 2}^\infty  {\frac{1}{{\nu \ln \nu }}}  =  + \infty .\]
This implies that  $  \mathcal{H}_{{{\log }^\alpha }}^{\tau+1}({\mathbb{T}^n})\not\subset\widetilde {\mathcal{L}_{\tau  + 1,b}^2}({\mathbb{T}^n})  $ for any $ 0<\alpha\leqslant1 $. This completes the proof of Theorem \ref{JZT2}.

\subsection{Proof of Theorem \ref{JZT3}}
For $ p = \frac{1}{2},p' = \frac{p}{{1 - p}} =  - 1 $, utilizing the reverse H\"older inequality (see \cite[Corollary 12]{MR3431607} for instance), we get
\begin{align*}
	\sum\limits_{{b^{\nu  - 1}} < \left| k \right| \leqslant {b^\nu }} {{{\left| {{f_k}} \right|}^2}{{\left| k \right|}^{2\tau  + 2}}}  &\geqslant {\left( {\sum\limits_{{b^{\nu  - 1}} < \left| k \right| \leqslant {b^\nu }} {{{\left| {{f_k}} \right|}^{2p}}} } \right)^{\frac{1}{p}}}{\left( {\sum\limits_{{b^{\nu  - 1}} < \left| k \right| \leqslant {b^\nu }} {{{\left| k \right|}^{\left( {2\tau  + 2} \right)p'}}} } \right)^{\frac{1}{{p'}}}}\\
	&= {\left( {\sum\limits_{{b^{\nu  - 1}} < \left| k \right| \leqslant {b^\nu }} {\left| {{f_k}} \right|} } \right)^2}{\left( {\sum\limits_{{b^{\nu  - 1}} < \left| k \right| \leqslant{b^\nu }} {\frac{1}{{{{\left| k \right|}^{2\tau  + 2}}}}} } \right)^{ - 1}}\\
	& \simeq {\left( {\sum\limits_{{b^{\nu  - 1}} < \left| k \right| \leqslant {b^\nu }} {\left| {{f_k}} \right|} } \right)^2}{\left( {\int_{{b^{\nu  - 1}}}^{{b^\nu }} {\frac{1}{{{r^{2\tau  + 3 - n}}}}{\rm d}r} } \right)^{ - 1}}\\
	& \simeq {\left( {\sum\limits_{{b^{\nu  - 1}} < \left| k \right| \leqslant{b^\nu }} {\left| {{f_k}} \right|} } \right)^2}{b^{\left( {2\tau  + 2 - n} \right)\nu }}.
\end{align*}
By summing up $ \nu $, we arrive at
\begin{align*}
	\sum\limits_{\nu  = 0}^\infty  {{{\left( {\sum\limits_{{b^{\nu  - 1}} < \left| k \right| \leqslant {b^\nu }} {{{\left| {{f_k}} \right|}^2}{{\left| k \right|}^{2\tau+2}}} } \right)}^{\frac{1}{2}}}} & \gtrsim \sum\limits_{\nu  = 0}^\infty  {\left( {\sum\limits_{{b^{\nu  - 1}} < \left| k \right| \leqslant{b^\nu }} {\left| {{f_k}} \right|} } \right){b^{\left( {\tau   - n/2+ 1} \right)\nu }}} \\
	& \simeq \sum\limits_{\nu  = 0}^\infty  {\left( {\sum\limits_{{b^{\nu  - 1}} < \left| k \right| \leqslant {b^\nu }} {\left| {{f_k}} \right|{{\left| k \right|}^{\left( {\tau   - n/2+ 1} \right)\nu }}} } \right)} \\
	& = \sum\limits_{k \in {\mathbb{Z}^n}} {\left| {{f_k}} \right|{{\left| k \right|}^{\left( {\tau   - n/2+ 1} \right)\nu }}} .
\end{align*}
This implies $ \widetilde {\mathcal{L}_{\tau  + 1,b}^2}({\mathbb{T}^n}) \subset {{C}^{\left[ {\tau   - n/2} \right]+ 1}}({\mathbb{T}^n}) $ in that for $ f \in  \widetilde {\mathcal{L}_{\tau  + 1,b}^2}({\mathbb{T}^n}) $,
\[\mathop {\sup }\limits_{x \in {\mathbb{T}^n}} \left| {{D^{\left[ {\tau  - n/2} \right]+ 1 }}f(x)} \right|  \lesssim  \sum\limits_{k \in {\mathbb{Z}^n}} {\left| {{f_k}} \right|{{\left| k \right|}^{\left[ {\tau   - n/2} \right]+ 1}}}  \leqslant \sum\limits_{k \in {\mathbb{Z}^n}} {\left| {{f_k}} \right|{{\left| k \right|}^{\left( {\tau   - n/2+ 1} \right) }}}  <  + \infty .\]
This completes the proof of Theorem \ref{JZT3}.

\subsection{Proof of Theorem \ref{JZT4} (Main Theorem I)}
The first claim follows directly from Theorem \ref{JZT3}, therefore it suffices to prove the latter.

We assert that $ \widetilde {\mathcal{L}_{\tau  + 1,b}^2}({\mathbb{T}^n}) \subset  \mathcal{H} ^{[2\tau]+2}({\mathbb{T}^n})  $. As for $ f(x)=\sum\nolimits_{k \in {\mathbb{Z}^n}} {{f_k}{{\rm e}^{{\rm i}\left\langle {k,x}  \right\rangle}}} \in \widetilde {\mathcal{L}_{\tau  + 1,b}^2}({\mathbb{T}^n})$, we have 
\[{\left( {\sum\limits_{{b^{\nu  - 1}} < \left| k \right| \leqslant {b^\nu }} {{{\left| {{f_k}} \right|}^2}{{\left| k \right|}^{2\tau+2}}} } \right)^{\frac{1}{2}}} \geqslant \left( {\sum\limits_{{b^{\nu  - 1}} < \left| k \right| \leqslant {b^\nu }} {{{\left| {{f_k}} \right|}^2}{{\left| k \right|}^{2\tau+2}}} } \right),\quad\text{as }\nu\gg 1,\]
hence, by summing up $ \nu $ and using $ 2\tau+2 \geqslant [2\tau]+2\geqslant 2n$, we prove the assertion.  This shows that every element in $ \widetilde {\mathcal{L}_{\tau  + 1,b}^2}({\mathbb{T}^n}) $ must possess $ n $-th order weak derivatives. In what follows, we construct an example in $  \widetilde {\mathcal{L}_{\tau  + 1,b}^2}({\mathbb{T}^n}) $ that admits unbounded weak derivatives from order $ \left[ {\tau    - n/2} \right]+2 $ to $  [\tau] +1 $. Construct $ \sum\nolimits_{k \in {\mathbb{N}^n}} {{f_k}\cos \left\langle {k,x} \right\rangle } $ with $ {f_k} = \left( {f_k^1, \ldots ,f_k^n} \right) $ with  $ f_k^s =f_k^t$ for all $ 1\leqslant s, t \leqslant n $, and let 
\[f_k^1 \sim \frac{1}{{{{\left| k \right|}^{\tau  + 1 + \frac{n}{2}}}{{\left( {\ln \left| k \right|} \right)}^\alpha }}},\quad\text{as }|k| \to  + \infty \]
for $ \alpha>1 $. Then  $ f(x) \in \widetilde {\mathcal{L}_{\tau  + 1,b}^2}({\mathbb{T}^n}) $ due to 
\[\sum\limits_{k \in {\mathbb{N}^n},|k| \geqslant 2} {{{\left| {{f_k}} \right|}^2}{{\left| k \right|}^{2\tau  + 2}}{{\left( {\ln \left| k \right|} \right)}^\alpha }}  \simeq \sum\limits_{k \in {\mathbb{N}^n},|k| \geqslant 2} {\frac{1}{{{{\left| k \right|}^n}{{\left( {\ln \left| k \right|} \right)}^\alpha }}}}  \simeq \int_2^{ + \infty } {\frac{{{r^{n - 1}}}}{{{r^n}{{\ln }^\alpha }r}}{\rm d}r}  = \int_2^{ + \infty } {\frac{1}{{r{{\ln }^\alpha }r}}{\rm d}r}  <  + \infty \]
and $  \mathcal{H}_{{{\log }^\alpha }}^{\tau+1}({\mathbb{T}^n})\subset\widetilde {\mathcal{L}_{\tau  + 1,b}^2}({\mathbb{T}^n})$ with $ \alpha>1 $ in Theorem \ref{JZT2}. We will demonstrate that this function satisfies the desired properties. Before proceeding, we establish an elementary inequality: for any  $ 0 \ne k \in {\mathbb{N}^n} $, it holds
\begin{equation}\label{CAU}
	\sum\limits_{j = 1}^n {k_j^\sigma }  \geqslant {n^{ - \sigma  - 1}}{\left| k \right|^\sigma }, \quad \sigma  \geqslant 1.
\end{equation}
Indeed, by utilizing the generalized Cauchy-Schwarz inequality, we obtain
\begin{align}
	\frac{{\sum\limits_{j = 1}^n {k_j^\sigma } }}{{{{\left| k \right|}^\sigma }}} &= \frac{{{{\left( {\frac{{{k_1}}}{{{k_1} +  \cdots  + {k_n}}}} \right)}^\sigma }}}{{{1^{\sigma  + 1}}}} +  \cdots  + \frac{{{{\left( {\frac{{{k_n}}}{{{k_1} +  \cdots  + {k_n}}}} \right)}^\sigma }}}{{{1^{\sigma  + 1}}}}\notag \\
	& \geqslant \frac{{{{\left( {\frac{{{k_1}}}{{{k_1} +  \cdots  + {k_n}}} +  \cdots  + \frac{{{k_n}}}{{{k_1} +  \cdots  + {k_n}}}} \right)}^\sigma }}}{{{{\left( {1 +  \cdots  + 1} \right)}^{\sigma  + 1}}}}\notag \\
\notag & = \frac{1}{{{n^{\sigma  + 1}}}}.
\end{align}
Now, applying \eqref{CAU} with $ \sigma={\left[ {\tau   - n/2} \right] + 2} \geqslant 1$, we arrive at 
\begin{align}
	\sum\limits_{j = 1}^n {\left( {\sum\limits_{k \in {\mathbb{N}^n}} {k_j^{\left[ {\tau   - n/2} \right] + 2}f_k^1} } \right)}  &\gtrsim \sum\limits_{k \in {\mathbb{N}^n},\;\left| k \right| \geqslant 2} {{{\left| k \right|}^{\left[ {\tau   - n/2} \right] + 2}}f_k^1} \notag \\
	& \simeq \sum\limits_{k \in {\mathbb{N}^n},\;\left| k \right| \geqslant 2} {\frac{1}{{{{\left| k \right|}^{\tau  + n/2 - \left[ {\tau   - n/2} \right]-1}}{{\left( {\ln \left| k \right|} \right)}^\alpha }}}} \notag \\
	& \simeq \int_2^\infty  {\frac{{{r^{n - 1}}}}{{{r^{\tau  + n/2 - \left[ {\tau  - n/2} \right]-1}}{{\ln }^\alpha }r}}{\rm d}r} \notag \\
	& = \int_2^\infty  {\frac{1}{{{r^{\tau   - n/2 - \left[ {\tau  - n/2} \right]}}{{\ln }^\alpha }r}}{\rm d}r} \notag \\
	\label{zuhe}& =  + \infty ,
\end{align}
where the last equality follows from the fact $ \tau    - n/2 - \left[ {\tau   - n/2} \right] \in [0,1) $. Finally, by using  \eqref{zuhe}, we obtain 
\[\mathop {\sup }\limits_{x \in {\mathbb{T}^n}} \left| {{D^{\left[ {\tau   - n/2} \right] + 2}}f(x)} \right| \gtrsim \sum\limits_{j = 1}^n {\left| {\partial_{{x_j}}^{\left[ {\tau    - n/2} \right] + 2}f\left( 0 \right)} \right|}  \gtrsim \sum\limits_{j = 1}^n {\left( {\sum\limits_{k \in {\mathbb{N}^n}} {k_j^{\left[ {\tau    - n/2} \right] + 2}f_k^1} } \right)}  =  + \infty ,\]
demonstrating  that the $ \left({\left[ {\tau    - n/2} \right] +2}\right) $-th order derivatives of $ f(x) $ can be unbounded (thus they are not classical derivatives, only of the weak type), and similarly for higher order (up to order $  [\tau] +1 $) weak derivatives.

It remains to demonstrate the existence of examples within $ \widetilde{\mathcal{L}_{\tau + 1, b}^2}(\mathbb{T}^n) $ that do not admit the $ ([ \tau ] + 2) $-th order weak derivatives. Consider constructing a suitable example \( f(x) = \sum_{0 \ne k \in \mathbb{Z}^n} f_k {\rm e}^{{\rm i} \langle k, x \rangle} \) that satisfies the condition:
\[
\sum_{b^{\nu - 1} < |k| \leqslant b^\nu} |f_k|^2 |k|^{2\tau + 2} := \Theta_\nu \sim \frac{1}{\nu^4}, \quad \text{as } \nu \to +\infty.
\]
It is evident that $ f(x) \in \widetilde{\mathcal{L}_{\tau + 1, b}^2}(\mathbb{T}^n) $, because
\[\sum\limits_{\nu  = 0}^\infty  {{{\left( {\sum\limits_{{b^{\nu  - 1}} < \left| k \right| \leqslant {b^\nu }} {{{\left| {{f_k}} \right|}^2}{{\left| k \right|}^{2\tau  + 2}}} } \right)}^{\frac{1}{2}}}}  = \sum\limits_{\nu  = 0}^\infty  {\Theta _\nu ^{\frac{1}{2}}}  \simeq \sum\limits_{\nu  = 1}^\infty  {\frac{1}{{{\nu ^2}}}}  <  + \infty .\]
However, \( f(x) \) does not admit the \( (\lfloor \tau \rfloor + 2) \)-th order weak derivatives due to the following observation:
\begin{align*}
	\sum\limits_{0 \ne k \in {\mathbb{Z}^n}} {{{\left| {{f_k}} \right|}^2}{{\left| k \right|}^{2\left( {\left[ \tau  \right] + 2} \right)}}} 	& = \sum\limits_{\nu  = 0}^\infty  {\sum\limits_{{b^{\nu  - 1}} < \left| k \right| \leqslant {b^\nu }} {{{\left| {{f_k}} \right|}^2}{{\left| k \right|}^{2\tau  + 2}} \cdot {{\left| k \right|}^{2\left( {\left[ \tau  \right] + 1 - \tau } \right)}}} } \\
	& \simeq \sum\limits_{\nu  = 0}^\infty  {{b^{2\left( {\left[ \tau  \right] + 1 - \tau } \right)\nu }}{\Theta _\nu }}  \\	&\simeq \sum\limits_{\nu  = 1}^\infty  {\frac{{{b^{2\left( {\left[ \tau  \right] + 1 - \tau } \right)\nu }}}}{{{\nu ^2}}}}  \\&=  + \infty .\quad \text{(note that $ \left[ \tau  \right] + 1 - \tau  > 0 $ for all $ \tau \geqslant 0 $)}
\end{align*}
Clearly, there are infinitely many such examples $ f(x) $. This completes the proof of Theorem \ref{JZT4}.

\subsection{Proof of Theorem \ref{TPDLINGLI4}}
For integral vectors $ \alpha  = \left( {{\alpha _1},{\alpha _2},\ldots,{\alpha _n}} \right) $ and $ k  = \left( {{k _1},{k _2},\ldots,{k _n}} \right) $, define $ {k^\alpha }: = k_1^{{\alpha _1}} \cdots k_n^{{\alpha _n}} $. Then for $ f\left( x \right) = \sum\nolimits_{0 \ne k \in {\mathbb{Z}^n}} {{f_k}{{\rm e}^{{\rm i}\left\langle {k,x} \right\rangle }}} \in 	\widetilde {\mathcal{L}_{\tau  + 1,b}^2}({\mathbb{T}^n}) $, we formally have  $ {\partial ^\alpha }f\left( x \right) = \sum\nolimits_{0 \ne k \in {\mathbb{Z}^n}} {{f_k}{{\rm i}^{\left| \alpha  \right|}}{k^\alpha }{{\rm e}^{{\rm i}\left\langle {k,x} \right\rangle }}}  $. Recall that we have proved in Theorem \ref{JZT4} that $ f(x) \in \mathcal{H}^{\tau+1} (\mathbb{T}^n) $.
It hence follows from  Parseval's
identity that 
\begin{equation}\label{OP111}
	\int_{{\mathbb{T}^n}} {{{\left| {{\partial ^\alpha }f\left( {y + s} \right) + {\partial ^\alpha }f\left( {y - s} \right) - 2{\partial ^\alpha }f\left( y \right)} \right|}^2}{\rm d}y}  = \sum\limits_{0 \ne k \in {\mathbb{Z}^n}} {4{{\left| {\cos \left\langle {k,s} \right\rangle  - 1} \right|}^2}{{\left| {{f_k}} \right|}^2}{{\left| {{k^\alpha }} \right|}^2}}.
\end{equation}
By setting $ s = \left( {l{b^{ - \nu - 1}},0} \right) $ with $ 0<l <\pi/2$  and $ k = ( {\bar k,\tilde k} ) $ with $ \bar k \in \mathbb{Z},\tilde k \in {\mathbb{Z}^{n - 1}} $, we obtain from \eqref{OP111} for all $ \nu \in \mathbb{N}^+ $:
\[\sum\limits_{0 \ne k \in {\mathbb{Z}^n},\;{b^\nu} < | \bar k | \leqslant {b^{\nu + 1}}} {{{\left| {{f_k}} \right|}^2}{{\left| {{k^\alpha }} \right|}^2}}  \simeq \sum\limits_{0 \ne k \in {\mathbb{Z}^n},\;{b^\nu} < | \bar k | \leqslant {b^{\nu + 1}}} {4{{\left| {\cos \left\langle {k,s} \right\rangle  - 1} \right|}^2}{{\left| {{f_k}} \right|}^2}{{\left| {{k^\alpha }} \right|}^2}}  \lesssim {\widetilde\varpi_{L^2} ^2}\left( {{\partial ^\alpha }f,{lb^{ - \nu-1}}} \right).\]
Proceeding with the same operations for the other components of $ k $ and summing up $ \alpha $ with $ |\alpha|=\tau+1$, we obtain
\begin{equation}\label{TPBDS2}
	\sum\limits_{\left| \alpha  \right| = \tau+1} {\sum\limits_{{b^{\nu-1}} < \left| k \right| \leqslant {b^{\nu  }}} {{{\left| {{f_k}} \right|}^2}{{\left| {{k^\alpha }} \right|}^2}} }  \lesssim \sum\limits_{\left| \alpha  \right| = \tau+1} {\widetilde\varpi _{{L^2}}^2\left( {{\partial ^\alpha }f,{b^{ - \nu}}} \right)}  = \widetilde\varpi _{{L^2}}^2\left( {{D^{\tau+1}}f,{lb^{ - \nu-1}}} \right).
\end{equation}
Note that by utilizing the generalized Cauchy-Schwarz inequality (similar to 
\eqref{CAU}), we obtain  
\begin{align}
	\sum\limits_{\left| \alpha  \right| = \tau+1} {\sum\limits_{{b^{\nu-1}} < \left| k \right| \leqslant {b^{\nu  }}} {{{\left| {{f_k}} \right|}^2}{{\left| {{k^\alpha }} \right|}^2}} }  &= \sum\limits_{{b^{\nu-1}} < \left| k \right| \leqslant {b^{\nu  }}} {{{\left| {{f_k}} \right|}^2}\sum\limits_{\left| \alpha  \right| = \tau+1} {{{\left| {{k^\alpha }} \right|}^2}} } \notag \\
	& \geqslant \sum\limits_{{b^{\nu-1}} < \left| k \right| \leqslant {b^{\nu  }}} {{{\left| {{f_k}} \right|}^2}\sum\limits_{\nu = 1}^n {\sum\limits_{{\alpha _\nu} = \tau+1} {{{\left| {{k^\alpha }} \right|}^2}} } }  \notag \\
	&= \sum\limits_{{b^{\nu-1}} < \left| k \right| \leqslant {b^{\nu  }}} {{{\left| {{f_k}} \right|}^2}\sum\limits_{\nu = 1}^n {\frac{{k_\nu^{2\tau+2}}}{{{1^{2\tau+1}}}}} }  \notag \\
	& \geqslant \sum\limits_{{b^{\nu-1}} < \left| k \right| \leqslant {b^{\nu  }}} {{{\left| {{f_k}} \right|}^2}{{\left( {\sum\limits_{\nu = 1}^n {|{k_\nu}|} } \right)}^{2\tau+2}}{{\left( {\sum\limits_{\nu = 1}^n 1 } \right)}^{ - \left( {2\tau+1} \right)}}}  \notag \\
	\label{TOBDS}& \simeq \sum\limits_{{b^{\nu-1}} < \left| k \right| \leqslant {b^{\nu  }}} {{{\left| {{f_k}} \right|}^2}{{\left| k \right|}^{2\tau + 2}}} .
\end{align}
Finally, combining \eqref{TPBDS2} and \eqref{TOBDS} yields
\begin{align}
	\sum\limits_{\nu = 0}^\infty  {{{\left( {\sum\limits_{{b^{\nu-1}} < \left| k \right| \leqslant {b^{\nu  }}} {{{\left| {{f_k}} \right|}^2}{{\left| k \right|}^{2\tau  + 2}}} } \right)}^{\frac{1}{2}}}}  &\lesssim \sum\limits_{\nu \geqslant 0} {{{\widetilde \varpi }_{{L^2}}}\left( {{D^{\tau+1}}f,l{b^{ - \nu - 1}}} \right)} 
	\notag \\
	& \simeq \sum\limits_{\nu \geqslant 0} {\left( {l{b^{ - \nu - 1}} - l{b^{ - \nu - 2}}} \right)\frac{{{{\widetilde \varpi }_{{L^2}}}\left( {{D^{\tau +1}}f,l{b^{ - \nu - 1}}} \right)}}{{l{b^{ - \nu - 1}}}}} \notag\\
\label{3434}	& \simeq \int_0^1 {\frac{{{{\widetilde \varpi }_{{L^2}}}\left( {{D^{\tau  + 1}}f,x} \right)}}{x}{\rm d}x} ,
\end{align}
where the last relation utilizes the monotonicity property of the $ (\tau+1)  $-th order Dini modified $ L^2 $-modulus of continuity, along with the following simple inequality:
\[\mathop {\sup }\limits_{x \in {\mathbb{T}^n}} \left| {{{\widetilde \varpi }_{{L^2}}}\left( {{D^{\tau  + 1}}f,x} \right)} \right| \leqslant 4{\left\| {{D^{\tau  + 1}}f} \right\|_{{L^2}\left( {{\mathbb{T}^n}} \right)}}.\]
This gives the proof of   Theorem \ref{TPDLINGLI4}.

\subsection{Proof of Theorem \ref{THPOSCHEL}  (Main Theorem II)} 
Recalling the finitely differentiable KAM theorem in \cite{POSCHELARXIV},  we deduce that the KAM conjugacy in \eqref{PTge} holds as long as the perturbation $ P(x) $ is small in the $ \| \cdot \|_{   \widetilde {\mathcal{L}_{\tau+1,b}^2}({\mathbb{T}^n})} $ norm. Following from Main Theorem I (Theorem \ref{JZT4}),  in this case, the perturbation $ P $ must possess $ (\left[ {\tau   - n/2} \right]+ 1) $-th order classical derivatives, but can admit unbounded weak derivatives from order $ \left[ {\tau   - n/2} \right]+2 $ to  $  [\tau]+1 $. Specifically, by setting $ \tau=n-1 $ we prove the first part of Theorem \ref{THPOSCHEL}.

 As for the second part, it  suffices to note that almost all vectors in $ \mathbb{R}^n $ have Diophantine exponent $ \tau $ for any $ \tau > n - 1 $. Let $ \epsilon = \tau - (n - 1)  $, and as $ \epsilon $ approaches $ 0 $ from the right, we have that $ \mathop {\lim }\nolimits_{\epsilon  \to {0^ + }} \left[ {n/2 + \epsilon } \right] = \left[ {n/2} \right] $ and $ \mathop {\lim }\nolimits_{\epsilon  \to {0^ + }} \left[ {n + \epsilon } \right] = n $. Thus, the second part of Theorem \ref{THPOSCHEL} is established.

To conclude the third part, it suffices to  combine the first part with Theorem \ref{TPDLINGLI4} for $ \tau :=n-1 \in \mathbb{N}^+ $ (see \eqref{3434} for details).

\subsection{Proof of Theorem \ref{sssss}}
In fact, the desired conclusion can be formulated in a neighborhood of $0\in \mathbb{R}^n$, since an example with  compact support on $\mathbb{R}^n$ can be shrunk down to ${\left[ { - 1/2,1/2} \right]^n}$ and periodized to produce an example on the torus. Vice versa, an example on the torus can be switched to an example on $\mathbb{R}^n$ by restricting on $ 2\pi $-period.
Without loss of generality, say $ x_0 =0 \in \mathbb{T}^n$. For any $ \sigma>1/2 $, construct a suitable example $ f(x) $ satisfying
\[{D^n}f\left( 0 \right) \simeq \frac{1}{{{{\left| x \right|}^{n/2}}{{\left( { - \ln \left| x \right|} \right)}^\sigma }}},\quad{\text{as }}\left| x \right| \to {0^ + }.\]
Here, we assume that $D^n f(x) $ only has a singularity point $ 0 $. Otherwise, for some
$C^\infty$ function $\varphi_n$ supported in $[-1,1]^n$ with \[\varphi_n \left( x \right) = \left\{ \begin{gathered}
	1, \quad 0 \leqslant \left| x \right| \leqslant 1/8, \hfill \\
	0, \quad 1/4 \leqslant \left| x \right| \leqslant 1/2, \hfill \\ 
\end{gathered}  \right.\]
it suffices to consider $ f^*(x) :=\varphi_n f(x) $ instead. 
With the choice of $ \sigma $, it is evident that $ f (x) \in \mathcal{H}^n (\mathbb{T}^n) $, because
\[\int_{\left| x \right| \leqslant \frac{1}{2}} {{{\left( {\frac{1}{{{{\left| x \right|}^{\frac{n}{2}}}{{\left( { - \ln \left| x \right|} \right)}^\sigma }}}} \right)}^2}{\rm d}x}  \simeq \int_0^{\frac{1}{2}} {\frac{{{r^{n - 1}}}}{{{r^{2n}}{{\left( { - \ln r} \right)}^{2\sigma }}}}{\rm d}r}  = \int_2^{ + \infty } {\frac{1}{{z{{\left( {\ln z} \right)}^{2\sigma }}}}{\rm d}z}  <  + \infty, \]
hence  $ f(x) \in \widetilde {\mathcal{L}_{ n,b}^2}({\mathbb{T}^n}) $ (see the proof of Theorem \ref{JZT4}). In this case, we have
\[{D^{n - 1}}f\left( 0 \right) \simeq \int_{\left| x \right|}^{\frac{1}{2}} {\frac{1}{{v_1^{\frac{n}{2}}{{\left( { - \ln {v_1}} \right)}^\sigma }}}{\rm d}{v_1}},\quad{\text{as }}\left| x \right| \to {0^ + }, \]
which proves \eqref{sss} with $ j=1 $. The other cases are similar.  This gives the proof of   Theorem \ref{sssss}.

 \section*{Acknowledgements} 
 Z. Tong extends heartfelt thanks to   Profs.  L. Evans,  L.  Grafakos and L. Slav\'{\i}kov\'{a} for valuable references on Sobolev spaces,  to Prof. L. Wang for valuable suggestions, and to Profs. A. Bounemoura  and M. Sevryuk for useful discussions.  Z. Tong  was supported by the China Postdoctoral Science Foundation (Grant No. 2025M783102). Y. Li was supported in part by the National Natural Science Foundation of China (Grant Nos. 12071175, 12471183 and 12531009).

\end{document}